\magnification=\magstep1
\input amstex
\documentstyle{amsppt}
\define\R{\Bbb R}
\define\rgen{\dot r_{gen}}
\pageheight{7in}

\topmatter
\title
Countable support iteration revisited
\endtitle
\author
Jind\v rich Zapletal
\endauthor
\affil
University of Florida
\endaffil
\abstract
The countable support iteration is the optimal way of iterating
proper definable forcings adding one real.
\endabstract
\thanks
The author is partially supported by grants GA \v CR 201-00-1466 and NSF DMS-0071437.
\endthanks
\address
Department of Mathematics, University of Florida, Gainesville FL 32611
\endaddress
\email
zapletal\@math.ufl.edu
\endemail
\subjclass
03E17, 03E55, 03E60
\endsubjclass
\endtopmatter

\document
\head {0. Introduction}\endhead

This paper is a sequel to \cite {Z}, using much of the same prerequisites and notation and proving a similar result. My purpose here is
to show that in a quite general setting, the countable support iteration is provably the best available tool for consistency results
concerning inequalities between cardinal invariants. The most frequent scenario for proving the consistency of such an inequality
$\frak x<\frak y$ is \cite {S} to first force the Continuum Hypothesis and then choose a proper forcing $P$ and iterate it with {\it countable support}
$\aleph_2$ many times, reaching the {\it iterated} $P$ {\it extension}. Presumably, the forcing $P$ was chosen in such a way that this extension
satisfies $\frak c=\frak y=\aleph_2$ and a {\it preservation theorem} was proved to guarantee that $\frak x=\aleph_1$ holds there.

The forcing $P$ is most frequently of the following syntactical kind \cite {B, RS}:

\definition {0.1. Definition}
$\langle P,\leq\rangle$ is a {\it definable real forcing} if there is a $P$-name $\dot r_{gen}$ for a real such that
\roster
\item"{$\circledast$}" the conditions in $P$ are real numbers, $P,\leq$ and $\dot r_{gen}$ are projective sets
\item"{$\circledast$}" for some projective formula $\chi(x,y)$ $P\Vdash$the generic filter is the set $\{p\in\check P:\chi(p,\dot r_{gen})\}$.
\endroster
\enddefinition

The equality $\frak y=\aleph_2$ in the iterated $P$ extension is most frequently obtained from the fact that $FA(P)$ holds in it, where

\definition {0.2. Definition}
For a partial order $P$ the forcing axiom for $P$, or $FA(P)$, is the statement that for every collection $\{O_\alpha:\alpha\in
\omega_1\}$ of $\aleph_1$ many open dense subsets
of $P$ and every condition $p\in P$ there is a filter $G\subset P$ containing the condition $p$ and meeting all the sets $O_\alpha,\alpha\in
\omega_1.$
\enddefinition

And the invariant $\frak x$ is usually tame as defined in \cite {Z}:

\definition{0.3. Definition}
A {\it tame} invariant is one defined as $\min\{|A|:A\subset\R, 
\phi(A)\land\psi(A)\}$ where $\phi(A)$ is a statement whose quatifiers range over the set $A$ and the natural numbers, and 
$\psi(A)$ is a statement of the form
$\forall x\in\R\ \exists y\in A\ \theta(x, y)$ where $\theta$ is a projective formula. A real parameter is allowed in both formulas $\phi$ and $\psi$.
\enddefinition

This is what I have to say about this general setup:

\proclaim {0.4. Theorem}
Suppose that there is a proper class of measurable Woodin cardinals. For every provably proper definable real forcing $P$ and every tame
invariant $\frak x$, if $FA(P)\land\frak x=\aleph_1$ holds in some forcing extension, then it holds in the iterated $P$ extension.
\endproclaim

Properly read, the theorem says that among a wide class of possible preservation properties, the countable support
iteration of a proper definable real forcing has every preservation property it can possibly have. This is certainly not a novel concept,
but it is interesting that it can be formally proved in this generality. Note also that $FA(P)$ can be obtained without iterating $P$ \cite {V}.

The assumption ``$FA(P)\land\frak x=\aleph_1$ holds in some forcing extension'' can be weakened to 
``in some extension, $\frak x$ is less than the smallest size of a family of open dense subsets of $P$ that cannot be met by a filter''.
For certain forcings $P$ (like $P=$Sacks) the assumption ``$FA(P)\land\frak x=\aleph_1$ in some extension''
can be weakened even further (to $\frak x<\frak c$ in some extension in the case of Sacks forcing, \cite {Z}).
However, a close study of that particular forcing $P$ is necessary for that.

The arguments can be retouched to yield at least one interesting variation of the theorem. If $P,Q$ are two proper definable real forcings
and $\frak x$ is a tame invariant such that $FA(P)\land FA(Q)\land\frak x=\aleph_1$ holds in some forcing extension then it holds
in the iterated $P*\dot Q$ extension. Roughly, if $FA(P)\land FA(Q)\land\frak x=\aleph_1$ is consistent then so is $FA(P*\dot Q)\land
\frak x=\aleph_1.$ This is curious inasmuch $FA(P)\land FA(Q)$ does not in general imply $FA(P*\dot Q).$

The arguments in this paper make use of the following three important unpublished results.

\proclaim {0.5. Fact}
(Woodin) ($\Sigma^2_1$ absoluteness) Suppose there is a proper class of measurable Woodin cardinals. Whenever $\phi$ is a $\Sigma^2_1$ sentence
with a possible real parameter which holds in some set forcing extension, then it holds in every set forcing extension satisfying
the continuum hypothesis.
\endproclaim

This fact is not explicitly used in this paper, but it does appear in \cite {Z, Section 2} to which this paper refers.

\proclaim {0.6. Fact}
(Woodin) (Real determinacy) Suppose that there is an inaccessible cardinal $\delta$ that is a limit of Woodin cardinals and $<\delta$-strong cardinals.
Then there is an inner model of ZF+DC+AD$\R$ containing all the reals such that all of its sets of reals are $<\delta$-weakly homogeneously
Souslin.
\endproclaim

\noindent This result is really an overkill for our purposes, however it provides a convenient way of sweeping certain definability issues under the rug
and of introducing the methodologically important models of real determinacy. All models as in Fact 0.6 will be called {\it good}. Their
internal theory is subject to 

\proclaim {0.7. Fact}
(Woodin) (Transfinite determinacy) In the base theory of ZF+DC,
 AD$\R$ is equivalent with the assertion that real games of any fixed countable length are determined.
\endproclaim

The following three fairly well known facts will be used in the paper without explicit notice.

\proclaim {0.8. Fact} 
\cite {W1} (Weakly homogeneous absoluteness)
 Suppose that $\delta$ is a limit of Woodin cardinals and $M$ is an elementary submodel of a large enough structure.
Then, in all generic extensions of $V$ via posets of size $<\delta$, all generic extensions of the model $M$ via posets of size $<\delta$
are correct about projective truth and the projections of $<\delta$-weakly homogeneous trees.
\endproclaim

\proclaim {0.9. Fact}
(Folklore? \cite {Z, Claim 1.1.3})
 Suppose that $P$ is a forcing adding a real and $M$ is a countable elementary submodel of a large enough structure. Then
the set $\{r\in\R:r$ is $M$-generic for the poset $P\cap M\}$ is Borel.
\endproclaim

\proclaim {0.10. Fact}
(Folklore?) Suppose that $J$ is a $\sigma$-ideal on the reals. Then the poset Borel($\R$) modulo $J$ adds a single real which is the
unique real belonging to all the sets in the generic filter.
\endproclaim

\noindent For the last fact, argue by transfinite induction on the Borel rank of the Borel positive set $X\subset\R$ that $X\Vdash\dot r\in\dot X$
where $\dot r$ is the real given by $\dot r(\check n)=\check m$ if the set $\{s\in\R: s(n)=m\}$ is in the generic filter.

Many lemmas in this paper are proved using this or that kind of a large cardinal hypothesis, determinacy hypothesis or a forcing absoluteness hypothesis.
After initial attempts at keeping exact account of these assumptions resulted in cumbersome mantra recited before every statement, I decided
to indicate their necessity by a simple (LC) in the wording of the lemmas. It is obvious from the proofs exactly what was used, however
the optimization of the hypotheses can be a difficult process. The assumptions of all results in this paper, except for one, could be
reduced to the existence of $\omega_1$ Woodin cardinals or less at the cost of introducing less intuitive notation and references to
several hard unpublished results of Martin, Neeman and Woodin. The only place where the measurable Woodin cardinals are used
is the $\Sigma^2_1$ absoluteness argument in \cite {Z, Subsection 2.3} that is referred to in the beginning of Section 3.

My notation follows closely the standard of \cite {J}. AD$\R$ is the statement that all games of length $\omega$ with real moves are determined.

\head {1. The single step forcing}\endhead

Fix a proper definable real forcing $P$ with its name $\dot r_{gen}$ for a generic real, and let $\chi(x, y)$ be a projective formula
such that $P\Vdash$ the generic filter $\dot G$ is the set $\{p\in\check P:\chi(p,\dot r_{gen})\}.$ Note that under $\omega$ Woodin
cardinals this implies that $V[r_{gen}]=V[G]$ since the model $V[r_{gen}]$ is correct about projective truth in $V[G]$ and therefore
computes the filter $G$ correctly.

 For a set $A\subset\R$ define a game $G(A)$ between
Adam and Eve in the following way: Eve starts, plays a condition $p_0\in P$ and then the two players alternate to produce a descending
chain $p_0\geq p_1\geq p_2\geq\dots$ in the poset $P.$ Eve wins if the filter $g$ generated by this chain meets all the antichains
necessary for the valuation of the name $\dot r_{gen}$, and $\dot r_{gen}/g\in A.$ Let $I_P$ be the collection of all set $A\subset\R$
for which Adam has a winning strategy in the game $G(A).$

\proclaim {1.1. Lemma}
$I_P$ is a $\sigma$-ideal.
\endproclaim

\demo {Proof}
The collection $I_P$ is clearly closed under subsets. Suppose now that $A_n,n\in\omega$ are sets in $I_P$ and $\sigma_n,n\in\omega$
are the relevant Adam's strategies. I must produce a winning strategy for Adam in the game $G(\bigcup_nA_n)$. Well, Adam will win
if he creates the following log along with his moves:

$$\matrix
p_0 & p_1 & p_2   & aux_0 & p_4   & aux_1 & p_6   & aux_3 & \dots\\
\   & \   & aux_0 & p_3   & aux_1 & aux_2 & aux_3 & aux_4 & \dots\\
\   & \   & \     & \     & aux_2 & p_5   & aux_4 & aux_5 &\dots\\
\   & \   & \     & \     & \     & \     & aux_5 & p_7   &\dots\\
\   & \   & \     & \     & \     & \     & \     & \     &\ddots
\endmatrix$$

Here the lines are simulated plays respecting the strategies $\sigma_0, \sigma_1,\dots$ so Eve's moves in these plays are recorded
in columns $0, 2, 4\dots$ The entries $p_0, p_1, p_2\dots$ are the ones actually made in the run of the game $G(\bigcup_nA_n)$ that
Adam is trying to win while the $aux$ entries are auxiliary. The only other rule for keeping the log is that
each simulated Eve's move is equal to the simulated Adam's move just up and right from it.

There is clearly exactly one way of keeping the log. The columns $2m$ and $2m+1$ are always obtained together from top to bottom, the
condition $p_{2m}$ is the weakest in these two columns and $p_{2m+1}$ is the strongest.
After the game is over, the filter $g$ generated by the descending chain of conditions actually played is equal to each of the filters
$g_n$ generated by the descending chain in the simulated play respecting the strategy $\sigma_n$ in the $n$-th row of the log.
Thus for each $n\in\omega,$ $\dot r_{gen}/g=\dot r_{gen}/g_n\notin A_n,$ so $\dot r_{gen}\notin\bigcup_nA_n$ and Adam has won. \qed
\enddemo
 
It is not hard to see that the ideal $I_P$ is invariant under forcing isomorphism. That is, if $P$ and $Q$ are forcings with their
respective generic reals $\dot r_{gen}$ and $\dot s_{gen}$ and $\pi$ is an isomorphism of their complete algebras carrying $\rgen$ to
$\dot s_{gen}$ then $I_P=I_Q.$

So far the definability of the forcing $P$ was not used. The main point in requiring the forcing to be definable is to make sure that many
of the above games are in fact determined. It will be instructive to see how the ideal $I_P$ behaves in inner models of real determinacy.

\proclaim {1.2. Lemma} (LC)
For every good model $N$ of ZF+DC+AD$\R$,
\roster
\item $(I_P)^N=I_P\cap N$
\item in $N,$ every set is either in $I_P$ or it has a Borel subset not in $I_P.$
\endroster
\endproclaim

\demo {Proof}
(1) is immediate. Since $N$ contains all the reals, the poset $P$ is in $N$ and if $N\models\sigma$ is a winning strategy for Adam or Eve
in a game $G(Y)$ then that strategy is winning even in $V.$ Since $N\models$AD$\R$, all the games $G(Y)$ for $Y\in N$ are determined with
relevant strategies in the model $N.$ So if $N\models Y\in I_P$ then $Y\in I_P$ because Adam's winning strategy in $N$ is still winning
in $V.$ And if $N\models Y\notin I_P$ then $N\models$ Eve has a winning strategy in $G(Y)$, this winning startegy is still winning in $V$ and
so certainly $Y\notin I_P.$

For (2) suppose that $Y\in N$ and $Y\notin I_P.$ Thus there is a suitably weakly homogeneous tree $T$ such that $Y=p[T]$ and there is a winning
strategy $\sigma\in N$ for Eve in the game $G(Y).$ Let $M$ be a countable elementary submodel of a large enough structure containing all the 
relevant information, and let $Z=\{r\in\R:r$ is an $M$-generic real for the poset $P$ and $M[r]\models r\in p[T]\}$. By Fact 0.9
the set $Z$ is Borel and since $Z\subset p[T],$ certainly $Z\subset Y.$ Thus to complete the proof it is enough to show that $Z\notin I_P.$
Let $\langle O_n:n\in\omega\rangle$ be an enumeration of all open dense subsets of $P$ in the model $M.$
Eve will in fact win the game $G(Z)$: as the play $p_0\geq p_1\geq p_2\geq\dots$ develops she keeps a log in the form $aux_0\geq
p_0\geq p_1\geq aux_1\geq aux_2\geq p_2\geq p_3\geq aux_3\geq\dots$ where the play consisting of the auxiliary moves follows the strategy $\sigma$
and the condition $aux_{2n+1}$ has an element of $O_n\cap M$ above it, and $p_0$ is any $M$-master condition below $aux_0=\sigma(0)\in M.$
This is easy to do, and in the end the filters $g$ and $h$ generated by the descending chain actually played and the chain of the auxiliary moves
will be the same. Since the strategy $\sigma$ was winning, $\rgen/g=\rgen/h\in Y.$ Moreover by the construction of the log, the real
$\rgen/g$ is $M$ generic for the poset $P,$ and since the tree $T$ was weakly homogeneous, $\rgen/g\in p[T]=Y$ implies $M[\rgen/g]\models\rgen/g\in p[T].$
Thus Eve won.

In fact (2) is a consequence of ZF+DC+AD$\R$. \qed
\enddemo

Now look at the poset $Q=$Borel($\R$) modulo the ideal $I_P.$ Thus Lemma 1.2(2) says that in any good model of real determinacy the poset $Q$
is dense in Power$(\R)$ modulo $I_P.$ By Fact 0.10, the poset $Q$ adds a single real $\rgen$ and every $I_P$-positive Borel set $X$
forces $\rgen$ into $\dot X.$ I can safely denote the generic real with $\rgen$ by the following representation theorem:

\proclaim {1.3. Lemma}
(LC) The poset $P$ is forcing isomorphic to $Q$ and the isomorphism takes the canonical $P$-generic real into the $Q$-generic real.
\endproclaim

\demo {Proof}
Let me first show that $Q\Vdash\dot G=\{p\in\check P:\chi(p,\rgen)\}$ is a generic filter on $P$. Suppose that $X$ is a Borel $I_P$ positive set,
$p_0$ and $p_1$ are conditions in $P$ such that $X\Vdash\check p_0,\check p_1\in\dot G$ and let $O\subset P$ be an open dense subset. I will produce
a set $Z\subset X$ and a condition $q\in P$ such that $q\in O$ and $q\leq p_0,p_1$, and $Z\Vdash\check q\in \dot G$. This will prove that $\dot G$
is forced to be a generic filter.

Let $M$ be a countable elementary submodel of a large enough structure containing all the relevant information and let $Y\subset X$ be the set
of all $M$-generic reals for the poset $P$ in the set $Y.$ Use the real projective determinacy to argue as in the previous proof that the set
$Y$ is $I_P$-positive; it is certainly Borel by Fact 0.9. By a projective forcing absoluteness argument $Y\Vdash\dot G\cap\check M\subset
\check P\cap\check M$ is an $\check M$-generic filter containing the conditions $\check p_0,\check p_1.$ So there must be a Borel $I_P$ positive
set $Z\subset Y$ and a condition $q\in O\cap M$ below both $p_0$ and $p_1$ such that $Z\Vdash\check q\in\dot G$ as desired.

To complete the proof of the lemma it is only necessary to show that for each condition $p\in P$ there is a Borel $I_P$ positive set
$X$ such that $X\Vdash\check p\in\dot G.$ Let $Y=\{r\in\R:\chi(p,r)\}.$ Clearly, $Y$ is a projective set and if it is $I_P$ positive then 
any Borel positive
subset $X$ of it will force $\check p$ into $\dot G$ by a projective absoluteness argument. Now Eve will in fact win in the game $G(Y)$:
as the play $p=p_0\geq p_1\geq p_2\geq\dots$ develops she will keep a log of countable elementary submodels $M_n, n\in\omega$ of some
large enough structure making sure that $M_0\in M_1\in\dots$, $p_n\in M_n$ and that in the end the filter $g$ generated by the descending
chain played is $M=\bigcup_nM_n$-generic. Obviously $p\in g,$ so $M[g]\models\chi(p,\rgen/g)$ and since $M[g]$ is correct
about projective truth, $\chi(p,\rgen/g)$ holds. Thus $\rgen/g\in Y$ and Eve won as desired. \qed
\enddemo

Thus it is safe to replace the forcing $P$ with $Q$. There are several points to be noted. First of all, $I_P=I_Q.$ However, the set of
(codes for) the conditions in $Q$ is no longer projective. Rather, it is in the definability class $\Game_{\R}$, but as such it will still belong
to any good model of ZF+DC+AD$\R$, and it will be weakly homogeneous provided suitable large cardinals exist. Of course the very forcing 
equivalence of $P$ and $Q$ was proved under a large cardinal hypothesis. The descriptive set theoretic complexity of the ideal $I_P$
is intimately connected to the forcing properties of $P$; several theorems exploring the relationship will appear elsewhere.

\head {2. The iteration reformulated}\endhead

The geometric representation of countable support countable length iterations of the forcing $P$ now proceeds as in \cite {Z, Section 1}.

\definition {2.1. Definition}
For an ordinal $\alpha\in\omega_1$ define the poset $P_\alpha$ to
consist of the nonempty Borel sets $p\subset\Bbb
R^\alpha$ satisfying these three conditions:
\roster
\item"{$\circledast$}" For every ordinal $\beta\in\alpha$ the set $p\restriction\beta=\{\vec s\in\Bbb
R^\beta:\exists \vec r\in p\ \vec s\subset\vec r\}$ is Borel. (The projection
condition)
\item"{$\circledast$}" For every ordinal $\beta\in\alpha$ and every sequence
$\vec s\in p\restriction\beta$, the set
$\{t\in\Bbb R:\vec s^\smallfrown\langle t\rangle\in p\restriction\beta+1\}$ is $I_P$-positive. (The $P$
condition)
\item"{$\circledast$}" For every increasing sequence $\beta_0\in\beta_1\in\dots$ of ordinals below $\alpha$
and every inclusion increasing sequence of sequences $\vec s_0\in p\restriction\beta_0, \vec
s_1\in p\restriction\beta_1\dots$, the sequence $\bigcup_n\vec s_n$ is in the set
$p\restriction\bigcup_n\beta_n$. (The countable support condition.)
\endroster
The sets $P_\alpha$ are ordered by inclusion.
\enddefinition 

Lemma 1.2 of \cite {Z} holds just the same, saying among other things that $P_\alpha$ is forcing equivalent to the countable support iteration
of the poset $P$ of length $\alpha.$ However, the lemma now requires a large cardinal assumption to guarantee throughout the argument
that internal generic extensions of countable elementary submodels evaluate the $I_P$-positivity correctly for Borel sets. If there is any good
model of ZF+DC+AD$\R$ then the set (of codes for conditions in) $Q$ is weakly homogeneous and the computation will be done correctly.
The weak homogeneity of $Q$ can be derived from the existence of $\omega+\omega$ Woodin cardinals.

The dichotomy \cite {Z, Lemma 1.4} takes a different form. Given a countable ordinal $\alpha$ and a set $A\subset\R^\alpha$ consider
the real game of length $\omega\alpha$ where

\roster
\item Adam and Eve play reals coding Borel sets $\langle Y_{\omega\beta+n}:\beta\in\alpha, n\in\omega\rangle$; at limit stages, Eve starts
\item for every ordinal $\beta\in\omega$ the sets $Y_{\omega\beta+n}:n\in\omega$ are $I_P$-positive and inclusion decreasing.
\endroster

Eve wins if for every ordinal $\beta\in\omega$ the set $\bigcap_nY_{\omega\beta+n}$ is a singleton $\{r_\beta\}$ and the sequence $\langle
r_\beta:\beta\in\alpha\rangle$ is in the set $A$. I will call the sequence $\langle
r_\beta:\beta\in\alpha\rangle$ the {\it result} of the play. I will also abuse the notation in assuming that Adam and Eve actually play the
Borel sets instead of reals coding them since except for the complexity calculations the difference will be immaterial,
and the resulting expressions will be shorter.

By the Transfinite Determinacy Fact 0.7, if there is any inner model of ZF + DC + AD$\R$ containing all the reals, 
the games $G(A)$ for projective sets $A\subset\R^\alpha$ must be
determined. We have

\proclaim {2.2. Lemma}
(LC) Suppose $\alpha\in\omega_1$ is a countable ordinal and $A\subset\R^\alpha$ is a projective set such that Eve has a winning
strategy in the game $G(A).$ Then there is a condition $p\in P_\alpha$ with $p\subset A.$
\endproclaim

\demo {Proof}
The proof follows the lines of the argument for \cite {Z, Claim 1.2.3}. Suppose that $A$ is a projective set such that Eve has a winning strategy
in the game $G(A)$, and suppose that there is a good model $N$ of ZF+DC+AD$\R$. By Fact 0.7, Eve must have a winning strategy $\sigma\in N.$
By transfinite induction on $\beta\leq\alpha$ argue that for any ordinal $\gamma\in\beta$, any $\gamma$-sequence $\vec r$
of reals and any partial play $\vec Y$ of the game $G(A)$ of length $\omega\gamma$ observing the strategy $\sigma$ whose result is $\vec r$
there is a condition $p\in P_{\beta-\gamma}$ and a function $f\in N$ with domain $p$ such that for every sequence $\vec s\in p$ the functional value
$f(\vec s)$ is a sequence such that $\vec Y^\smallfrown f(\vec s)$ is a partial play of the game $G(A)$ of length $\omega\beta$
observing the strategy $\sigma$, whose result is the sequence $\vec r^\smallfrown\vec s.$ This will prove the lemma considering the case
$\beta=\alpha,\gamma=0$ and the fact that $\sigma$ is a winning strategy for Eve.

Let me just perform the successor step of the induction. Suppose $\beta=\beta'+1$ and the induction hypothesis has been verified up to $\beta'$,
suppose $\gamma\in\beta,$ $\vec r\in\R^\gamma$ and $\vec Y$ is a partial play of the game $G(A)$ observing the strategy 
$\sigma$ resulting in the sequence $\vec r.$ Find a condition $p'\in P_{\beta'-\gamma}$ and a function $f'\in N$
as in the induction hypothesis for $\beta'$, for $\gamma=\beta'$ this will be $p'=f'=0.$ For every sequence $\vec s\in p'$ let $C_{\vec s}=\{t\in\R:$ 
for some sequence $\vec X$ the sequence $\vec Y^\smallfrown f'(\vec s)^\smallfrown\vec X$ is a partial play
of length $\omega\beta$ respecting the strategy $\sigma$ resulting in the sequence $\vec r^\smallfrown\vec s^\smallfrown\langle t\rangle\}.$
The set $C_{\vec s}$ is $I_P$-positive: this follows from the fact that $I_P=I_Q$ and that a suitable fraction of the strategy
$\sigma$ gives Eve a winning strategy in the game $G(C_{\vec s}).$ I will find a condition $p\in P^{\beta-\gamma}$ which is a subset
of the set $\{\vec s^\smallfrown\langle t\rangle:\vec s\in p'\land t\in C_{\vec s}\}$. This will complete the induction step because then
Uniformization in $N$ can be used to find a function $f\in N$ with domain $p$ such that $f(\vec s^\smallfrown\langle t\rangle)=f'(\vec s)^\smallfrown
\vec X$ where $\vec X$ is a witness to $t\in C_{\vec s}$, and this is exactly what is required in the induction hypothesis.

To find the set $p,$ note that the set $\{\langle \vec s, c\rangle:\vec s\in p'$ and $c$ is a code for a Borel $I_P$-positive subset
of the set $C_{\vec s}\}$ is in the model $N,$ and as such it is the projection of a suitably weakly homogeneous tree $T.$ 
By Weakly Homogeneous Absoluteness Fact 0.8,
$p'\Vdash_{P_{\beta'-\gamma}}$ for some code $c$ $\langle\vec r_{gen},c\rangle\in p[\check T]$; fix a name $\dot c$ for such a code.
Let $M$ be a countable elementary submodel of a large enough structure containing all the relevant information and find a condition
$q'\subset p'$ in the poset $P_{\beta'-\gamma}$ consisting solely of $M$-generic sequences of reals. Finally, let
$p=\{\vec s^\smallfrown\langle t\rangle\in\R^{\beta-\gamma}:\vec s\in q'$ and $t$ belongs to the Borel set coded by $\dot c/\vec s\}.$ \qed
\enddemo

\head {3. The absoluteness argument}\endhead

The rest of the proof of Theorem 0.4 follows closely the lines of Section 2 in \cite {Z}. I will only prove the following lemma which will replace the
argument for (*) in Subsection 2.1 of \cite {Z}. Let $\frak x$ be an arbitrary tame invariant defined as $\min\{|A|:A\subset\R, 
\phi(A)\land\psi(A)\}$ where $\phi(A)$ quantifies over $A$ and $\omega$ only and 
$\psi(A)$ is a statement of the form
$\forall x\in\R\ \exists y\in A\ \theta(x, y)$ where $\theta$ is a projective formula.

\proclaim {3.1. Lemma}
(LC) Suppose $FA(P)\land\frak x=\aleph_1$ holds and let $A\subset\R$ be any set of size $\aleph_1$ with the property $\psi(A)$.
Then for every ordinal $\alpha\in\omega_1,$ every condition $p\in P_\alpha$ and every Borel function $\dot x:p\to\R$
there is a condition $q\leq p$ and a real $y\in A$ such that for every sequence $\vec r\in q,$ $\theta(\dot x(\vec r),y)$ holds.
\endproclaim

To prove the lemma, choose an ordinal $\alpha\in\omega_1,$ a condition $p\in P_\alpha$ and a Borel function $\dot x:p\to\R.$
For each real $y\in A$ let $B_y=\{\vec r\in p:\theta(\dot x(\vec r), y)\}.$ These sets are projective, therefore the games $G(B_y)$
take place in any good model of ZF+DC+AD$\R$ and as such are determined. There are two cases.

If Eve has a winning strategy in one of the games $G(B_y)$ for a real $y\in A$ then by Lemma 2.2 there is a condition $q\subset B_y.$
Obviously, the real $y\in A$ and the condition $q\leq p$ are the desired objects.

Thus the proof will be complete if I derive a contradiction from the assumption that Adam has a winning strategy $\sigma_y$ in every game
$G(B_y), y\in A.$ By simultaneous transfinite induction on $\beta\leq\alpha$ build plays $\{ Y_{\omega\gamma+n,y}:\gamma\in\beta, n\in\omega,
y\in A\}$ of the games
$G(B_y)$ respecting the strategies $\sigma_y$ so that for each $\gamma\in\beta$ the unique real contained in the intersection
$\bigcap_nY_{\omega\gamma+n,y}$ does not depend on $y,$ and denoting this real by $r_\gamma,$ for each ordinal $\beta\leq\alpha$ the sequence 
$\langle r_\gamma:\gamma\in\beta\rangle$ belongs to the set $p\restriction\beta.$ Once this construction is complete, look
at the sequence $\vec r=\langle r_\beta:\beta\in\alpha\rangle.$ Since the strategies $\sigma_y$ were all winning for Adam
and we have just produced plays that resulted in the sequence $\vec r$ for every real $y\in A,$ it must be the case that
$\vec r\notin B_y$ for all $y\in A.$ However, this means that for all $y\in A$ $\theta(\dot x(\vec r),y)$ fails, contradicting
the property $\psi$ of the set $A$.

To realize the transfinite induction mentioned in the previous paragraph, suppose the plays $\{Y_{\omega\gamma+n,y}:\gamma\in\beta, n\in\omega, y\in A\}$
have been constructed for an ordinal $\beta\in\alpha.$ For every real $y\in A$ by induction on $n\in\omega$ build maximal
antichains $E_{yn}\subset Q$ and maps $f_{yn}$ with domain $E_{yn}$ so that

\roster
\item $E_{yn+1}$ refines $E_{yn}$ and for every set $Y\in E_{yn}$, the value $r(n)$ is the same for all reals $r\in Y$
\item for every set $Y\in E_{yn}$ the functional value $f_{yn}(Y)$ is some 
decreasing sequence $Y_0\geq Y_1\geq Y_2\geq \dots \geq Y_{2n+1}=Y$
of sets in the forcing $Q$ such that the play $\langle Y_{\omega\gamma+m,y}:\gamma\in\beta,m\in\omega\rangle^\smallfrown\langle Y_m:m\in 2m+2
\rangle$ follows the strategy $\sigma_y$
\item whenever $Z\subset Y$ are sets in the respective antichains $E_{yn}$ and $E_{ym}$ where $m\in n$, necessarily $f_{ym}(Y)
\subset f_{yn}(Z).$
\endroster

The induction step of this construction is handled in the following way. Once $E_{yn}$ and $f_{yn}$ have been constructed, for each
set $Y\in E_{yn}$ consider the set $D_Y=\{ Z\subset Y:$ for some set $X\subset Y$ the play $\langle Y_{\omega\gamma+m,y}:\gamma\in\beta,
m\in\omega\rangle^\smallfrown f_{yn}(Y)^\smallfrown\langle X, Z\rangle$ follows the strategy $\sigma_y$, and moreover, the value
$r(n+1)$ is the same for every real $r\in X\}$. By the definition of the game $G(B_y)$ this set is dense in $Q$ below $Y.$
Now just let $E_{yn+1}$ be any maximal antichain of $Q$ included in the union $\bigcup_{Y\in E_{yn}}D_Y$ and for each set 
$Z$ in it, let $Y$ be the unique superset of $Z$ in $E_{yn},$ let $X\subset Y$ be the witness for $Z\in D_Y$ and let
$f_{yn+1}(Z)=f_{yn}(Y)^\smallfrown\langle X, Z\rangle.$

Now the collection $\{ E_{yn}:n\in\omega, y\in A\}$ of maximal antichains of $Q$ has size $\aleph_1.$ Use $FA(P)=FA(Q)$ and \cite {SZ, Lemma 38}
to find an elementary submodel $M$ of size $\aleph_1$ of a large enough structure containing all the relevant information
and this collection and the set $A$ as subsets, and an $M$-generic filter $G_\beta\subset Q\cap M$ containing the condition
$\{r\in\R:\langle r_\gamma:\gamma\in\beta\rangle^\smallfrown\langle r\rangle\in p\restriction\beta+1\}\in Q$. 
By Fact 0.10 applied in the model $M$ the intersection
$\bigcap G_\beta$ contains exactly one real $r_\beta.$ For each real $y\in A$ and each number $n\in\omega$ there is a unique element
$Z_{yn}$ in the intersection $E_{yn}\cap G_\beta$. By (1-3) above, the union $\bigcup_n f_{yn}(Z_{yn})$ is an $\omega$-sequence
$\langle Y_{\omega\beta+n,y}:n\in\omega\rangle$ whose intersection is the singleton $\{r_\beta\}$ and the partial play
$\langle Y_{\omega\gamma+m,y}:\gamma\in\beta,m\in\omega\rangle^\smallfrown \langle Y_{\omega\beta+n,y}:n\in\omega\rangle$
follows the strategy $\sigma_y.$ This completes the inductive step in the simultaneous transfinite induction on $\beta\in\alpha,$ and the proof
of the lemma.

\enddocument